# Symbolic estimation of distances between eigenvalues of Hermitian operator and unitary orbits classification

Ilia Lomidze, Natela Chachava (Georgian Technical University, Tbilisi, Georgia)
E-mail: ilia.lomidze@gtu.ge, chachava.natela@yahoo.com

**Abstract.** We use symbolic expressions for traces of positive integer powers of a Hermitian operator (or, equivalently, coefficients of corresponding characteristic polynomial) to find solutions for the problems as follows:

Factorization of characteristic polynomial that collects all eigenvalues that have the same multiplicity into one polynomial factor with coefficients rationally expressed through coefficients of original characteristic polynomial.

Finding with any accuracy a minimal distance between eigenvalues, as well as a minimal and a maximal eigenvalues of the Hermitian operator, applying only rational functions of corresponding coefficients of the characteristic polynomial. Estimation of related rate of convergence.

Classification of unitary orbits of Hermitian operators.

**Mathematics Subject Classification** (2010). 12Y99, 47B15.

## 1. Introduction.

Our goal is to describe the class of unitary equivalent Hermitian matrices by their unitary invariants, namely, by traces of their natural degrees. All matrices that are related to each other by a similarity transformation

$$P = UQU^{-1} \quad \Leftrightarrow \quad P \sim Q,$$

where $U$ is a unitary $N \times N$ matrix, obviously, have the same set of these invariants and form a class of equivalency – the orbit of the correspondent Hermitian operator. The necessary condition for unitary equivalence (that is whether the matrices $P$ and $Q$ belong to the same orbit) is obvious: the matrices $P$ and $Q$ must have the same spectrum, that is, the same set of eigenvalues taking into account their multiplicities. For Hermitian matrices this condition is also sufficient. Hence, two problems arise, as follows:

1) Considering a certain set of unitary invariants to find out whether the matrices $P$ and $Q$ belong to the same orbit;

2) To classify various orbits, that is, to describe equivalence classes, using a certain set of unitary invariants.

The problem of checking whether two or more given Hermitian matrices belong to the same orbit was solved in [1], where the following results were obtained: two given Hermitian matrices $P$ and $Q$ belong to the same orbit if and only if the following traces coincide:

$$\mathrm{tr} P^k = \mathrm{tr} Q^k = \sum_{i=1}^{m} r_i p_i^k, \ k = \overline{0, 2m-1}.$$

Here $m$ is the number of different eigenvalues of the matrix (the number of different roots of the corresponding characteristic polynomial, CP). It was shown in [1] that between these $2m$ conditions there exist $m - l$ syzygies arising if there are eigenvalues with the same multiplicity; $l$ stands for the number of different multiplicities. The issues above, related to the traditional theory of matrix invariants, became useful in some applications. In particular, in classification of states of quantum computer [2] by means of the invariants of corresponding (Hermitian) density matrix the following problem occurred: how to split the set of orbits into equivalence classes and describe these classes by unitary invariants of the matrix. In present paper we propose splitting of the orbits set into classes in the following way: orbits of two Hermitian (generally, nonsimilar) matrices

$$P \sim \operatorname{diag}\{\underbrace{p_1,\ldots,p_1}_{r_1},\ldots,\underbrace{p_m,\ldots,p_m}_{r_m}\},\ p_1 < \ldots < p_m$$

and

$$Q \sim \operatorname{diag}\{\underbrace{q_1,\ldots,q_1}_{r_1},\ldots,\underbrace{q_m,\ldots,q_m}_{r_m}\},\ q_1 < \ldots < q_m$$

belong to the same class. Hence, we label classes of orbits by a set $\{r_1,\ldots,r_m\}$ that corresponds to the ordered sequence of eigenvalues of these matrices. Actually, the question is: what conditions (syzygies) are satisfied by unitary invariants of matrices belonging to given class.

While studying the problem we have obtained even more detailed description of Hermitian matrix spectrum: we have factorized CP of the matrix collecting all eigenvalues having the same multiplicities in one factor – polynomial with coefficients rationally expressed through CP coefficients. Besides, we have found (with any accuracy) minimal distance between eigenvalues, as well as minimal and maximal eigenvalues of the operator, using only rational functions of CP coefficients. Then we have splitted the region containing all eigenvalues of CP onto uniform cells of maximal length such that each cell contains at most one eigenvalue. We can check whether the certain cell is empty or not using only rational functions of CP coefficient. Hence, we can classify orbits in accordance with the numbers of occupied/empty cells using only rational functions of CP coefficient.

## 2. Notations

Let us consider polynomial $P_n(x) \in \mathbb{R}[x]$ having real roots only – $p_1 < p_2 < \ldots < p_m$ – with corresponding multiplicities $r_1, r_2, \ldots, r_m$ ($r_i \in \mathbb{N}$, $\sum_{i=1}^{m} r_i = n$). In what follows we use notations as in the article [1]):

$$t_k = \sum_{l=1}^{m} r_l p_l^k,\ k = 0, 1, 2, \ldots, \qquad (1)$$

$$H_k(P_n) = H_k = \left[t_{i+j-2}\right]_1^k, \qquad (2)$$

$H_k$ is $k$-order Hankel matrix [3]-[5]. Taking into account that

$$t_{i+j-2} = \sum_{l=1}^{m} p_l^{i-1} r_l p_l^{j-1} = \left[p_l^{i-1}\right]^T \left[r_l p_l^{j-1}\right],$$

($\left[p_l^{i-1}\right]$ denotes one-column matrix, and $\left[p_l^{i-1}\right]^T$ – correspondent one-row matrix) and using the Binet-Cauchy Theorem (see [5]) it is easy to check (see [1]) that

$$D_k = \det H_k = \begin{cases} \sum_{1 \le i_1 < \cdots < i_k \le m} \left\{ r_{i_1} \cdots r_{i_k} \prod_{1 \le j < l \le k} (p_{i_j} - p_{i_l})^2 \right\} > 0, & 1 \le k \le m, \\ 0, & k > m. \end{cases} \qquad (3)$$

In particular, one has

$$D_2 = \det H_2 = \det \begin{bmatrix} t_0 & t_1 \\ t_1 & t_2 \end{bmatrix} = \det\left(\begin{bmatrix} 1 & \ldots & 1 \\ p_1 & \ldots & p_m \end{bmatrix}\begin{bmatrix} r_1 & \ldots & r_m \\ r_1 p_1 & \ldots & r_m p_m \end{bmatrix}^T\right) = \sum_{1 \le i < j \le m} r_i r_j (p_j - p_i)^2, \qquad (4)$$

$$D_m = \det H_m = \left(\prod_{i=1}^{m} r_i\right) \prod_{1 \le i < j \le m} (p_j - p_i)^2. \qquad (5)$$

Note, that all sums (1) and therefore all determinants (3)-(5) can be rationally expressed via coefficients of the polynomial, using known Newton's formulae. The explicit formulae (3)-(5) illustrate the statement (see [4]-[6]) of general

**Theorem** (Jacoby): *A number of different roots* (ndr) *of a polynomial* $P_n(x) \in \mathbb{R}[x]$ *is equal to the rank of the Hankel matrix* $H_n$ *(see (2))*:

$$\mathrm{ndr} = \mathrm{rank}(H_n);$$

*a number of different real roots* (ndrr) *of a polynomial* $P_n(x) \in \mathbb{R}[x]$ *is equal to the Hankel matrix signature*

$$\mathrm{ndrr} = P(1, D_1, ..., D_n) - V(1, D_1, ..., D_n).$$

Here a symbol $P(1, D_1, ..., D_n)$ (a symbol $V(1, D_1, ..., D_n)$) denotes a number of sign permanencies (a number of sign variations) in the sequence $1, D_1, ..., D_n$. Obviously, for the polynomial having real roots only, according to the inequalities (3), one has $V(1, D_1, ..., D_n) = 0$, $\mathrm{rank}(H_n) = P(1, D_1, ..., D_n) = m$ and $\mathrm{ndr} = \mathrm{ndrr} = m$.

## 3. Multiplicities

In the article [1] it was shown that any $m$ equalities from (1) can be solved for $r_i$, $i = \overline{1, m}$:

$$r_i^{-1} = D_m^{-1} \det \begin{bmatrix} H_m & (\mathbf{p}_i)^T \\ -\mathbf{p}_i & 0 \end{bmatrix}, \tag{6}$$

where

$$\mathbf{p}_i \equiv \left[ p_i^{k-1} \right]_{k=\overline{1,m}} = \left[ 1, p_i, ..., p_i^{m-1} \right]$$

Stands for a row, whereas $\mathbf{p}_i^T$ denotes a correspondent column.

Calculating the determinant in (6), we obtain

$$r_i^{-1} = \left\langle \mathbf{p}_i H_m^{-1} \mathbf{p}_i^T \right\rangle, \tag{7}$$

where $H_m^{-1}$ denotes the inverse matrix and brackets stand for a scalar product of a row by a column.

**Statement 1**. The formula (7) can be generalized as:

$$r_i^{-1} \delta_{ij} = \left\langle \mathbf{p}_i H_m^{-1} \mathbf{p}_j^T \right\rangle, \tag{8}$$

where $\delta_{ij}$ denotes the Kronecker delta.

**Proof.** It is enough to present the inverse matrix $H_m^{-1}$ as

$$H_m^{-1} = \left( \left[ p_i^{i-1} \right] \left[ r_l p_l^{j-1} \right]^T \right)^{-1} = \left( \left[ r_l p_l^{j-1} \right]^{-1} \right)^T \left[ p_i^{i-1} \right]^{-1}$$

and then calculate the scalar product (8).

**Corollary 1**. The next formula is valid

$$r_j^{-1} - r_i^{-1} = \left\langle (\mathbf{p}_j + \mathbf{p}_i) H_m^{-1} (\mathbf{p}_j - \mathbf{p}_i)^T \right\rangle. \tag{9}$$

**Proof.** Obvious.

Expanding the difference of degrees in the last column in the formula (9) we obtain

$$p_j^{k-1} - p_i^{k-1} = (p_j - p_i) \sum_{l=0}^{k-2} p_i^{k-2-l} p_j^l, \quad k = \overline{2, m}. \tag{10}$$

Let us construct the expression

$$\frac{r_j^{-1} - r_i^{-1}}{p_j - p_i} \equiv \Delta_{ij} = \Delta_{ji} = \sum_{n,k=1}^{m} (p_j^{n-1} + p_i^{n-1}) [H_m^{-1}]_{nk} \sum_{l=0}^{k-2} p_i^{k-2-l} p_j^l \tag{11}$$

(In the right hand side of the last formula the summand with $k=1$ is equal to zero). Sum of all $\Delta_{ij}$ ($i, j$, $1 \leq i < j \leq m$), obviously, gives us the symmetric function of $p_k$ in the right hand side of the (11).

The symmetric function can be expressed via the elementary symmetric functions of the variable $p_k$

$$\sigma_1 = \sum_{1 \leq j \leq m} p_j, \quad \sigma_2 = \sum_{1 \leq j < k \leq m} p_j p_k, \ldots, \quad \sigma_m = p_1 p_2 \cdots p_m, \quad (12)$$

(or, alternatively, via sums $s_k = \sum_{i=1}^{m} p_i^{k-1}$, $k = \overline{1, m}$, using the Newton's formulae).

For forthcoming consideration we need the next theorem:

**Theorem 1.** *All elementary symmetric functions (12) can be rationally expressed via sums (1) and therefore via the coefficients of the given polynomial $P_n(x)$.*

**Proof.** Let us consider the polynomial

$$P_m(x) = D_m^{-1} \det \begin{bmatrix} t_0 & t_1 & \cdots & t_{m-1} & t_m \\ t_1 & t_2 & \cdots & t_m & t_{m+1} \\ \cdots & \cdots & \cdots & \cdots \\ t_{m-1} & t_m & \cdots & t_{2m-2} & t_{2m-1} \\ 1 & x & \cdots & x^{m-1} & x^m \end{bmatrix} = D_m^{-1} \det \begin{bmatrix} H_m & (\mathbf{t})^T \\ \mathbf{x} & x^m \end{bmatrix}, \quad (13)$$

having the first coefficient 1. Let us show that other coefficients of this polynomial are from the set (12).

Let us use the obvious decomposition

$$\begin{bmatrix} H_m & (\mathbf{t})^T \\ \mathbf{x} & x^m \end{bmatrix} = \begin{bmatrix} [p_l^{i-1}]_1^m & \mathbf{0}^T \\ \mathbf{0} & 1 \end{bmatrix} \begin{bmatrix} ([r_l p_l^{j-1}]_1^m)^T & [r_l p_l^m]^T \\ \mathbf{x} & x^m \end{bmatrix}.$$

Then, taking into account the formulae (3)-(5) and the Viet theorem one easily finds

$$P_m(x) = D_m^{-1} \det [p_l^{i-1}]_1^m \left(\prod_{l=1}^{m} r_l\right) \det \begin{bmatrix} ([p_l^{j-1}]_1^m)^T & [p_l^m]^T \\ \mathbf{x} & x^m \end{bmatrix} =$$

$$= D_m^{-1} \left(\prod_{l=1}^{m} r_l\right) \left(\det [p_l^{k-1}]\right)^2 \prod_{i=1}^{m}(x - p_i) = \prod_{i=1}^{m}(x - p_i) =$$

$$= x^m + \sum_{k=1}^{m} x^{m-k}(-1)^k \sigma_k.$$

Hence, one gets

$$\sigma_k = \sum_{1 \leq i < j \leq m} p_{i_1} \cdots p_{i_k} = D_m^{-1} H\binom{m+1}{m-k+1}, \quad (14)$$

where $H\binom{m+1}{m-k+1}$ denotes the minor that remains after deleting the $(m+1)$-th row and the $(m-k+1)$-th column, $k = \overline{1, m}$, of the determinant in the formula (13). These minors, as well as the factor $D_m^{-1}$, are polynomial functions of the sums (1) and therefore the statement of the Theorem 1 follows. ∎

**Corollary 2**. The polynomial (13) is a minimal polynomial of Hermitian operator (matrix) having the polynomial $P_n(x)$ as a characteristic one.

**Corollary 3**. The function

$$Z \equiv \sum_{i<j} \frac{r_j^{-1} - r_i^{-1}}{p_j - p_i} = \sum_{i<j} \left\{ \sum_{n,k=1}^{m} (p_j^{n-1} + p_i^{n-1})[H_m^{-1}]_{nk} \sum_{l=0}^{k-2} p_i^{k-2-l} p_j^{l} \right\} \quad (15)$$

is a rational function of the coefficients of the polynomial $P_n(x)$.

**Remark 1**. The formula (13) can be rewritten in the form

$$\det \begin{bmatrix} H_m & (\mathbf{t})^T \\ \mathbf{x} & x^m \end{bmatrix} = D_m P_m(x) = \mathcal{H}_m(x), \quad (16)$$

that corresponds to the known result (see [4]) in a special case $r_i = 1$, $i = \overline{1, n}$, $m = n$; in contrast with the formula referred our result takes into account multiplicities of roots explicitly.

## 4. Distance between roots

Let us apply the formula (15) to the polynomial

$$Q_{3m}(x) = P_m(x) P_m^2(x - \varepsilon) = \prod_{i=1}^{m}(x - p_i) \prod_{j=1}^{m}(x - (p_j + \varepsilon))^2, \quad (17)$$

assuming that the $\varepsilon$ is chosen according to the condition

$$p_i - (p_j + \varepsilon) \neq 0, \quad i, j = \overline{1, m}, \, m \geq 3. \quad (18)$$

Applying the definition (15) to the polynomial (17) one gets

$$Z(\varepsilon) = \left(1 - \frac{1}{2}\right) \sum_{i,j=1}^{m} \frac{1}{p_i - (p_j + \varepsilon)} = \frac{1}{2} \left( \sum_{i<j} \frac{1}{p_i - (p_j + \varepsilon)} - \frac{m}{\varepsilon} + \sum_{i>j} \frac{1}{p_i - (p_j + \varepsilon)} \right) =$$

$$= \frac{1}{2} \left( -\frac{m}{\varepsilon} + \sum_{1 \leq i < j \leq m} \frac{2\varepsilon}{(p_i - p_j)^2 - \varepsilon^2} \right). \quad (m \geq 3) \quad (19)$$

Let us assume that

$$\varepsilon = \varepsilon_0, \quad 0 \leq \varepsilon_0 < \min |p_i - p_j| \equiv \mu.$$

Then one obtains

$$\frac{1}{\varepsilon_0} \left( Z(\varepsilon_0) + \frac{m}{2\varepsilon_0} \right) = \sum_{1 \leq i < j \leq m} \frac{1}{(p_i - p_j)^2 - \varepsilon_0^2} > \frac{1}{\mu^2 - \varepsilon_0^2} > 0. \quad (20)$$

As result, the next estimation follows to the inequality (20):

$$\mu^2 > \left( \sum_{1 \leq i < j \leq m} \frac{1}{(p_i - p_j)^2 - \varepsilon_0^2} \right)^{-1} + \varepsilon_0^2 \equiv \varepsilon_1^2. \quad (21)$$

Obviously,

$$\varepsilon_0^2 < \varepsilon_1^2 < \mu^2.$$

Hence, one can use $\varepsilon_1$ instead $\varepsilon_0$ in the formula (19) to calculate $\varepsilon_2 > 0$:

$$\varepsilon_2^2 \equiv \left( \sum_{1 \leq i < j \leq m} \frac{1}{(p_i - p_j)^2 - \varepsilon_1^2} \right)^{-1} + \varepsilon_1^2 < \mu^2.$$

For similar reasons one gets

$$\varepsilon_0^2 < \varepsilon_1^2 < \varepsilon_2^2 < \mu^2.$$

Continuing this process, after $k$ steps one obtains the increasing and bounded above sequence

$$\varepsilon_0^2 < \varepsilon_1^2 < \varepsilon_2^2 < ... < \varepsilon_k^2 < \mu^2, \tag{22}$$

where

$$\varepsilon_k^2 \equiv \left( \sum_{1 \le i < j \le m} \frac{1}{(p_i - p_j)^2 - \varepsilon_{k-1}^2} \right)^{-1} + \varepsilon_{k-1}^2, \quad k = 1, 2, .... \tag{23}$$

Therefore, when $k \to \infty$ the sequence (22) must converge to some limit $\lim_{k \to \infty} \varepsilon_k = \varepsilon_\infty$, $0 < \varepsilon_\infty \le \mu$. This limit can be easily calculated using the recurrent relation (23): assuming $k \to \infty$ in the both sides of it, one obtains

$$\varepsilon_\infty^2 \equiv \left( \sum_{1 \le i < j \le m} \frac{1}{(p_i - p_j)^2 - \varepsilon_\infty^2} \right)^{-1} + \varepsilon_\infty^2.$$

Hence,

$$\left( \sum_{1 \le i < j \le m} \frac{1}{(p_i - p_j)^2 - \varepsilon_\infty^2} \right)^{-1} = 0, \qquad \sum_{1 \le i < j \le m} \frac{1}{(p_i - p_j)^2 - \varepsilon_\infty^2} \to \infty.$$

So, in the last sum at least one summand is infinite – at least one denominator in it is equal to 0. Taking into account that $0 < \varepsilon_\infty \le \mu \equiv \min|p_i - p_j|$, one has to conclude

$$\varepsilon_\infty = \mu = \min|p_i - p_j|. \tag{24}$$

In order to perform this scheme really one has to prove the next

**Lemma 1.** *The parameters $\varepsilon_k$, $0 \le \varepsilon_k < \mu$, $k = 0, 1, 2, ...$, can be expressed through known coefficients of the polynomial $P_n(x)$.*

**Proof.** One can start with initial $\varepsilon_0 = 0$. Then, according to the definition (21) one gets

$$\varepsilon_1^2 = \left[ \sum_{1 \le i < j \le m} (p_i - p_j)^{-2} \right]^{-1} < \mu^2. \tag{25}$$

Obviously, the sum in the formula (25) is a symmetric function of the roots $p_1, ..., p_m$ of the minimal polynomial (13). Therefore, according to the Theorem 1, it can be expressed as a rational function of the elementary symmetric functions (14) and, hence, as a rational function of the coefficients of the given polynomial $P_n(x)$. Hence, the statement of the Lemma is valid for $\varepsilon_1$. Then, according to the recurrence formula (23) the Lemma can be proved using the mathematical induction for any $\varepsilon_k$, $k = 2, 3, ....$  ∎

So, we have proven the next statement:

**Statement 2.** *The minimal distance between the different roots of the polynomial $P_n(x)$ (having real roots only) can be find out as a limit of an increasing convergent sequence of the rational functions of the polynomial coefficients only.*  ∎

**Remark 2.** According to the Lemma 1, all terms of the convergent sequences $\{\varepsilon_k | k = 0, 1, ...\}$ are symmetric functions of the roots $p_i$, $i = \overline{1, m}$, $m \ge 3$, of the polynomial. It seems interesting to investigate whether the limit of the sequence $\varepsilon_\infty = \mu = \min|p_i - p_j|$ is also symmetric function of the roots $p_i$, $i = \overline{1, m}$, or not.

**Remark 3**. The function $Z(\varepsilon)$ defined by the formula (19) presents in the subsequent consideration only as a combination

$$\varepsilon^{-1}\left( Z(\varepsilon) + m(2\varepsilon)^{-1} \right) = \sum_{1 \le i < j \le m} \left[ (p_i - p_j)^2 - \varepsilon^2 \right]^{-1}.$$

So, one needs to express by coefficients of given polynomial the last symmetric function only. This is much easier than to express the function $Z(\varepsilon)$ based on the formulas (17) and (15).

Let us consider some $l$- subset $\omega_l$

$$\omega_l = \{i_1,...,i_l \mid l \leq m\} \subset \{1,...,m \mid m \leq n\}$$

and let $f(x_1,...,x_l)$ be a real function of $l \leq m$ real variables, bounded in the domain $x_j \in [p_1, p_m]$, $j=\overline{1,l}$:

$$\min\{f(x_1,...,x_l) \mid x_j \in [p_1, p_m], \ j=\overline{1,l}\} = a, \quad \max\{f(x_1,...,x_l) \mid x_j \in [p_1, p_m], \ j=\overline{1,l}\} = b.$$

One can prove the next

**Lemma 2**. *The recurrence relation*

$$\left[\sum_{\omega_l \in \{\overline{1,m}\}} \frac{1}{f(p_{i_1},...,p_{i_l}) - c_k}\right]^{-1} + c_k = c_{k+1}, \quad k=0,1,..., \tag{26}$$

*where $k = 1,2,...$, and real $c_0$ is assumed to be some known constant, gives us a sequence $\{c_k \mid k \in \mathbb{N}\}$ which has the next properties*:

a. *if $c_0 < a = \min\{f(x_1,...,x_l) \mid x_j \in [p_1, p_m], \ j=\overline{1,l}\}$, then the sequence is increasing and bounded above*:

$$c_0 < c_1 < ... < c_k < c_{k+1} < a; \tag{27}$$

b. *if $c_0 > b = \max\{f(x_1,...,x_l) \mid x_j \in [p_1, p_m], \ j=\overline{1,l}\}$, then the sequence is decreasing and bounded bellow*:

$$c_0 > c_1 > ... > c_k > c_{k+1} > b. \tag{28}$$

c. *All reals $c_k$, $k = 1,2,...$, defined by the recurrence relation (26) can be expressed through known $c_0$ and the coefficients of the polynomial $P_n(x)$.*

**Proof.** The statements *a)* and *b)* can be proved directly. To prove *c)* let us take into account that the expression

$$\sum_{\omega_l \in \{\overline{1,m}\}} [f(p_{i_1},...,p_{i_l}) - c_0]^{-1}$$

is a symmetric function of the $P_m(x)$ polynomial's roots. Then, according to the Theorem 1, one concludes that the real $c_1$, that is defined as

$$c_1 = \left[\sum_{\omega_l \in \{\overline{1,m}\}} \frac{1}{f(p_{i_1},...,p_{i_l}) - c_0}\right]^{-1} + c_0,$$

is a rational function of the coefficients of the polynomial $P_n(x)$ and $c_0$. Then, using the mathematical induction one can prove the statement *c)* of the Lemma 2 for any $c_k$, $k = 2,3,...$. ∎

**Corollary**. Assuming that $c_0$ is expressed by a rational function of the coefficients of the polynomial $P_m(x)$ one concludes that all terms $\{c_k \mid k \in \mathbb{N}\}$ of sequences (27) and (28) are rational functions of the coefficients of the polynomial $P_n(x)$.

**Remark 4**. Results obtained above for the minimal distance between the roots of the polynomial $P_m(x)$ (see formulae (22)-(24)) can be derived from the Lemma 2, choosing in (26) $f = f(x_1, x_2) = (x_1 - x_2)^2$ and $c_0 = 0$.

Let us consider the Hankel matrix $H_2(P_m) = [s_{i+j-2}]_{i,j=1,2}$, where

$$S_k = \sum_{i=1}^{m} p_i^{k-1}, \quad k = \overline{1,m},$$

are Newton's sums of the polynomial $P_m(x)$ roots. The determinant of this matrix

$$D_2(P_m) = s_0 s_2 - s_1^2 = m(\sigma_1^2 - 2\sigma_2) = \sum_{1 \le i < j \le m} (p_j - p_i)^2$$

can be obtained from the formula (5) putting there $r_i = 1, i = \overline{1,m}$. Obviously,

$$D_2(P_m) = \sum_{1 \le i < j \le m} (p_j - p_i)^2 > \max\{(p_j - p_i)^2 \mid i, j = \overline{1,m}\}.$$

Let us denote $\overline{p}(P_m)$ the arithmetic mean of the roots of the polynomial $P_m(x)$:

$$\overline{p}(P_m) = s_1 / m = m^{-1} \sum_{i=1}^{m} p_i.$$

**Lemma 3.** *The next estimations are valid:*

$$\min\{p_i \mid i = \overline{1,m}\} > \overline{p}(P_m) - \sqrt{D_2(P_m)} \qquad (29)$$

$$\max\{p_i \mid i = \overline{1,m}\} < \overline{p}(P_m) + \sqrt{D_2(P_m)}. \qquad (30)$$

**Proof.** Obvious. ∎

Let us choose in the formula (26) $f = f(x) = x$. Then one has

$$a = \min\{f(x) \mid x \in [p_1, p_m]\} = p_1,$$
$$b = \max\{f(x) \mid x \in [p_1, p_m]\} = p_m.$$

Now, according to the Lemma 2, the minimal root $p_1$ (the maximal root $p_m$) can be find out as a limit of increasing and bounded above sequence (27) (decreasing and bounded bellow sequence (28)), if one chooses $c_0$ according to

$$c_0 = \overline{p}(P_m) - \sqrt{D_2(P_m)}$$

(correspondingly, according to $c_0 = \overline{p}(P_m) + \sqrt{D_2(P_m)}$). Let us denote terms of the corresponding sequences by

$$c_0 = a_0 < a_1 < \ldots < a_k < a_{k+1} < a \qquad (31)$$
$$(c_0 = b_0 > b_1 > \ldots > b_k > b_{k+1} > b). \qquad (32)$$

Here

$$a_{k+1} = \left( \sum_{1 \le i \le m} \frac{1}{p_i - a_k} \right)^{-1} + a_k,$$

$$b_{k+1} = \left( \sum_{1 \le i \le m} \frac{1}{p_i - b_k} \right)^{-1} + b_k.$$

**Statement 4.** *Any term of the sequences (31) and (32) is symmetric function of the roots $p_i$, $i = \overline{1,m}$, and therefore is expressible as a rational function of the polynomial coefficients.*

**Proof.** Obvious. ∎

**Statement 5.** *All roots $p_1 < \ldots < p_m$ of the polynomial $P_m(x)$ (all different roots of the polynomial $P_n(x)$) are located in the open interval $]a_i, b_k[$, where $a_i$ and $b_k$ are some terms from the sequences (31) and (32) correspondingly.*

**Proof.** Obvious. ∎

The Statement 5 in some sense supplements the next (see [4], [7])

**Theorem** (Joachimsthal). *Let* $\text{rank}(H_n) = m$ *and let the sequence*

$$1, \mathcal{H}_1(x), ..., \mathcal{H}_m(x)$$

*does not contain two consecutive zeros for* $x = a$ *and* $x = b$. *Then*

$$\text{nrr}\{P_n(x) = 0 \mid a < x < b\} = V(1, \mathcal{H}_1(a), ..., \mathcal{H}_m(a)) - V(1, \mathcal{H}_1(b), ..., \mathcal{H}_m(b)).$$

(Hankel polynomials $\mathcal{H}_k(x)$, $k = \overline{1, m}$, are defined in (16)).

Namely, the Statement 5 specifies the smallest interval containing *all* roots of the polynomial $P_m(x)$ (all different roots of the CP $P_n(x)$).

## 5. Convergence rate

Let us estimate convergence rates of the sequence $\{\varepsilon_k \mid k = 0, 1, ...\}$.

Taking into account that for any $k = 0, 1, ...$ the inequalities $\varepsilon_k < \mu$, are valid, one gets estimations:

$$\mu^2 - \varepsilon_k^2 > \varepsilon_{k+1}^2 - \varepsilon_k^2 = \left( \sum_{1 \leq i < j \leq m} \frac{1}{(p_i - p_j)^2 - \varepsilon_k^2} \right)^{-1} \geq \left( \sum_{l=1}^{m-1} \frac{m-l}{l^2 \mu^2 - \varepsilon_k^2} \right)^{-1}. \quad (33)$$

The sums in the right hand side of the formula (33) correspond to the polynomial

$$W_m(x) = \prod_{l=0}^{m-1} (x - p_1 - \mu l), \quad \mu > 0, \quad (34)$$

having equidistant roots with the minimal distance between them $\mu$. Let us call this polynomial "Wilkinson's generalized polynomial" (WGP). The Wilkinson's polynomial [8] corresponds to the case $\mu = 1$, $p_1 = 1$, $p_m = m$. According to the estimation (33), WGP (34) provides the lowest rate of convergence of all sequences under consideration. For WGP a lot of estimations can be performed explicitly.

**Statement 6**. *For WGP (34) and for any* $k = 0, 1, ...$ *one has*

$$\varepsilon_k^2 = \mu^2 w_k^2(m), \quad (m \geq 3) \quad (35)$$

*and the coefficients* $w_k^2(m)$ *satisfy the recurrence relation*

$$w_{k+1}^2(m) = w_k^2(m) + \left( \sum_{l=1}^{m-1} \frac{m-l}{l^2 - w_k^2(m)} \right)^{-1}. \quad (36)$$

**Proof.** By induction. ∎

**Statement 7**. *One has*

$$0 = w_0^2 < w_1^2 < ... < w_k^2 < 1, \quad (37)$$

$$\lim_{k \to \infty} w_k(m) = 1, \quad (38)$$

$$w_k^2(m+1) < w_k^2(m), \quad (m \geq 3; \; k = 1, 2, ...) \quad (39)$$

**Proof.** (37) and (38) are obvious; (39) can be proved by induction. ∎

For WGP one has an equality in the right hand side of the formula (33) and hence one gets

$$w_{k+1}^2(m) - w_k^2(m) = \left\{ \sum_{l=1}^{m-1} \frac{m-l}{l^2 - w_k^2} \right\}^{-1} = \left\{ \frac{m-1}{1 - w_k^2} + \sum_{l=2}^{m-1} \frac{m-l}{l^2 - w_k^2} \right\}^{-1}.$$

Let us denote

$$1 - w_k^2 = v_k, \quad k = 0, 1, \ldots, \quad v_0 = 1, \quad v_k \underset{k \to \infty}{\to} 0+, \qquad (40)$$

and rewrite the previous formula as

$$v_k(m) - v_{k+1}(m) =$$

$$\left\{ \frac{m-1}{v_k} + \sum_{l=2}^{m-1} \frac{m-l}{l^2 - 1 + v_k} \right\}^{-1} = \frac{v_k}{m-1} \left\{ 1 + v_k \sum_{l=2}^{m-1} \frac{m-l}{m-1} \frac{1}{l^2 - 1} \left( 1 + \frac{v_k}{l^2 - 1} \right)^{-1} \right\}^{-1},$$

$$(v_0 = 1; \ k = 1, 2, \ldots). \qquad (41)$$

**Statement 8.**

$$\left( \frac{m-2}{m-1} \right)^k < v_k(m) \le \left[ 1 - A(m) \right]^k, \qquad (42)$$

where we denote

$$A(m) = \frac{(m-1)}{(m-1)^2 + m(m-2)B(m)} < 1,$$

$$B(m) = \sum_{l=2}^{m-1} \frac{m-l}{m-1} \frac{1}{l^2 - 1} = \frac{1}{m-1} \left( \frac{3}{4} m - \frac{1}{4} + \frac{1}{2m} - \sum_{l=1}^{m} l^{-1} \right) =$$

$$= \frac{3}{4} - \frac{1}{m-1} \left( 1 + \frac{1}{2m} + \sum_{l=3}^{m-1} l^{-1} \right). \qquad (m \ge 3) \qquad (43)$$

**Proof.** Obvious ∎

**Corollary.** Assuming that the accuracy precision of calculations is $\delta$, the convergence rate in (40) can be estimated based on (42) and on the inequality

$$v_k < \delta \ll 1$$

as

$$\frac{\ln \delta}{\ln(m-2) - \ln(m-1)} < k < \frac{\ln \delta}{\ln[1 - A(m)]}.$$

Besides, one has

$$B(3) = \frac{1}{6}; \quad B(4) = \frac{19}{72} \left( \approx \frac{1}{4} \right); \quad B(5) = \frac{79}{240} \left( \approx \frac{1}{3} \right); \quad B(m) \underset{m \to \infty}{\to} \frac{3}{4} - 0.$$

So, for $m \gg 1$ one gets

$$\ln(m-2) - \ln(m-1) \simeq -\frac{1}{m-2},$$

$$\ln[1 - A(m)] \simeq -A(m) \simeq -\frac{1}{m-1} \frac{1}{1 + \frac{m(m-2)}{(m-1)^2} \frac{3}{4}} \simeq -\frac{4}{7(m-1)},$$

and therefore the estimation gets the form

$$(m-2) \ln \delta^{-1} < k < \frac{7}{4} (m-1) \ln \delta^{-1}.$$

**Exsamples.**

$$m = 3, \quad A(3) = \frac{2}{4 + 3B(3)} = \frac{4}{9}, \quad 1 - A(3) = \frac{5}{9}; \quad 2^{-k} < v_k(3) \le (5/9)^k$$

$$\Rightarrow \quad k \ln(9/5) < \ln \delta^{-1} < k \ln 2; \quad \frac{\ln \delta^{-1}}{\ln 2} < k < \frac{\ln \delta^{-1}}{\ln(9/5)}.$$

$$m = 4, \quad A(4) = \frac{3}{9 + 8B(4)} = \frac{27}{100}, \quad 1 - A(4) = \frac{73}{100}, \quad \left(\frac{2}{3}\right)^k < v_k(4) < \left(\frac{73}{100}\right)^k.$$

$$\Rightarrow \quad k \ln(100/73) < \ln \delta^{-1} < k \ln(3/2); \quad \frac{\ln \delta^{-1}}{\ln(3/2)} < k < \frac{\ln \delta^{-1}}{\ln(100/73)}.$$

$$m \gg 1, \quad \delta \simeq \exp(-10) \quad \Rightarrow \quad 10(m-2) < k < 18(m-1).$$

## 6. Factorisation

Let us consider CP of a matrix

$$P \sim \mathrm{diag}\{\underbrace{p_1, ..., p_1}_{r_1}, ..., \underbrace{p_m, ..., p_m}_{r_m}\}, \quad p_1 < ... < p_m,$$

and let us assume that

$$r_{j_1} = ... = r_{j_{n_1}} \equiv q_1, ...,$$

$$r_{j_{m-n_l+1}} = ... = r_{j_m} \equiv q_l;$$

$$q_\alpha \neq q_\beta, \quad \alpha \neq \beta; \quad \sum_{\alpha=1}^{l} n_\alpha = m.$$

**Statement 9.** *The following factorization is fulfilled:*

$$P_n(x) = \prod_{\alpha=1}^{l} \left[Q_{n_\alpha}(x)\right]^{q_\alpha}, \tag{44}$$

*where*

$$Q_{n_\alpha}(x) = \prod_{\alpha=1}^{n_\alpha}(x - p^{\{q_\alpha\}}_i) = x^{n_\alpha} + \sum_{k=1}^{n_\alpha} x^{n_\alpha - k}(-1)^k \sigma^{\{q_\alpha\}}_k, \tag{45}$$

*and all coefficients $\sigma^{\{q_\alpha\}}_k$ can be rationally expressed by coefficients of CP $P_n(x)$. The factor $Q_{n_\alpha}(x), \alpha = \overline{1, l}$, has roots $p^{\{q_\alpha\}}_i, i = \overline{1, n_\alpha}$, all having the same multiplicity $q_\alpha$ in original CP $P_n(x)$.*

**Proof.** Let us consider the following power summands for some $q < n$ ($q \in \mathbb{N}$):

$$t_k^{(q)} = \sum_{i=1}^{m}(r_i - q)p_i^k = t_k - qs_k = \sum_{\alpha=1}^{l}(q_\alpha - q)\sum_{i=1}^{n_\alpha}(p^{\{q_\alpha\}}_i)^k. \tag{46}$$

As far all $t_k, s_k$ can be expressed by coefficients of the CP $P_n(x)$, we conclude that all $t_k^{(q)}$ are expressible by coefficients of the $P_n(x)$ too. If corresponding Hankel matrix

$$H_m^{(q)} = \left[t_{i+j-2}^{(q)}\right]_1^m$$

Is nonsingular, that is, if $m_q \equiv \mathrm{rank}(H_m^{(q)}) = m$, then there are no roots with multiplicity $q$; if for some $q_1$ one gets $m_{q_1} < m$, then CP $P_n(x)$ has $m - m_{q_1} \equiv n_1 < m$ different roots having multiplicity $q_1$. Performing this scenario for $q = 1, 2, ...$, one can find the minimal $q = q_1$ and use power summands (46) to construct the corresponding auxiliary polynomial (AP) $P_{m q_1}^{(q_1)}(x) = P_{m-n_1}^{(q_1)}(x)$:

$$P_{m q_1}^{(q_1)}(x) = D_{m-n_1}^{-1} \det \begin{bmatrix} t_0^{(q_1)} & t_1^{(q_1)} & \cdots & t_{m-n_1}^{(q_1)} \\ \cdots & \cdots & \cdots & \cdots \\ t_{m-n_1-1}^{(q_1)} & t_{m-n_1}^{(q_1)} & \cdots & t_{2(m-n_1)-1}^{(q_1)} \\ 1 & x & \cdots & x^{m-n_1} \end{bmatrix} = x^{m-n_1} + \sum_{k=1}^{m-n_1} x^{m-n_1-k}(-1)^k \sigma_k^{(q_1)}, \tag{47}$$

having all roots of the polynomial $P_m(x)$ (see (13)) except $p^{\{q_1\}}_k$, $k = \overline{1, n_1}$. Therefore, the first factor in (44) can be found as

$$Q_{n_1}(x) = \frac{P_m(x)}{P^{(q_1)}_{m-n_1}(x)} = \prod_{i=1}^{n_1}(x - p^{\{q_1\}}_i). \tag{48}$$

To find $Q_{n_2}(x)$ one has to construct

$$t_k^{(q_1, q_2)} = t_k^{(q_1)} - q_2 s_k^{(q_1)}, \quad (q_2 > q_1) \tag{49}$$

using AP (47) to construct $s_k^{(q_1)}$. Repeating the previous scenario one gets AP $P^{(q_1, q_2)}_{m-(n_1+n_2)}(x)$ and

$$Q_{n_2}(x) = \frac{P^{(q_1)}_{m-n_1}(x)}{P^{(q_1, q_2)}_{m-(n_1+n_2)}(x)} = \prod_{i=1}^{n_2}(x - p^{\{q_2\}}_i). \tag{50}$$

Continuing the process above one gets the partitions of integers $m$ and $n$ as follows:

$$m = \sum_{\alpha=1}^{l} n_\alpha, \quad n = \sum_{\alpha=1}^{l} n_\alpha q_\alpha. \tag{51}$$

Statement 9 now becomes obvious applying induction procedure with respect to $\alpha$. ∎

## 7. Syzygies

As far

$$\text{rank} H_m^{(q_\alpha)} = \text{rank}\left[t_{i+j-2}^{(q_\alpha)}\right]_1^m = m - n_\alpha$$

for any $q_\alpha$, $\alpha = \overline{1, l}$, one has the $n_\alpha$ equations

$$\det H_{m-n_\alpha+1}^{(q_\alpha)} = \ldots = \det H_m^{(q_\alpha)} = 0. \tag{52}$$

Excluding $q_\alpha$ from equations (52) one obtains $n_\alpha - 1$ syzygies. Total number of such syzygies, in agreement with result obtained in [1], is

$$\sum_{\alpha=1}^{l}(n_\alpha - 1) = m - l. \tag{53}$$

## 8. Lattice

Let us choose some term $\varepsilon_k \equiv \delta$ from (22) and some terms $a_i$, $b_k$ from the sequences (31), (32):

$$\delta < \mu = \min|p_i - p_j|,$$

$$a_i < p_1 = \min\{p_j \mid P_n(p_j) = 0\}, \quad b_k > p_m = \max\{p_j \mid P_n(p_j) = 0\}.$$

Let us define

$$x_j = a_i + j\delta, \quad j = \overline{0, M}, \tag{54}$$

where

$$M = [(b_k - a_i)/\delta + 1] \geq m,$$

($[x]$ denotes the integer part of a $x \in \mathbb{R}$). Obviously,

$$x_0 = a_i < p_1, \quad x_M = a_i + M\delta = a_i + [(b_k - a_i)/\delta + 1]\delta > b_k > p_m,$$

and the union of pairwise disjoint intervals

$$[x_0, x_1] \bigcup_{j=\overline{2, M}} \,]x_{j-1}, x_j]$$

includes the interval $[p_1, p_m]$. Besides, each cell $]x_{j-1}, x_j]$, $j = \overline{1, M}$, contains at most one root of the polynomial $P_m(x) = \prod_{\alpha=1}^{l} Q_{n_\alpha}(x)$. Therefore, calculating values of the polynomials at the sites of obtained lattice $Q_{n_\alpha}(x_j)$, $j = \overline{0, M}$, $\alpha = \overline{1, l}$, one can easily check the condition

$$Q_{n_\alpha}(x_{j-1}) Q_{n_\alpha}(x_j) \leq 0, \tag{55}$$

That is whether the given cell contains a root of the polynomial $Q_{n_\alpha}(x_j)$, for each $\alpha = \overline{1, l}$ separately. Note, that if $Q_{n_\alpha}(x_j) = 0$ for some $j = j_\beta$, $\beta = 1, 2, ...$, then this cell contains no other roots.

The set of natural numbers labeling nonempty cells $]x_{j-1}, x_j]$

$$\{j \mid Q_{n_\alpha}(x_{j-1}) Q_{n_\alpha}(x_j) \leq 0\} \equiv \{j(Q_{n_\alpha}(x); q_\alpha)\}, \tag{56}$$

depends on $\alpha = \overline{1, l}$, that is, depends on $q_\alpha$. Therefore, comparing the sets (56) for different $\alpha = \overline{1, l}$, one can determine the mutual arrangement of CP's roots (eigenvalues) having certain multiplicities $q_\alpha$.

**Definition.** *Orbits of two Hermitian (generally, nonsimilar) matrices*

$$P \sim \mathrm{diag}\{\underbrace{p_1, ..., p_1}_{r_1}, ..., \underbrace{p_m, ..., p_m}_{r_m}\}, \quad p_1 < ... < p_m$$

*and*

$$Q \sim \mathrm{diag}\{\underbrace{q_1, ..., q_1}_{r_1}, ..., \underbrace{q_m, ..., q_m}_{r_m}\}, \quad q_1 < ... < q_m$$

*are considered to belong to the same class, labeled by a set $\{r_1, ..., r_m\}$ that corresponds to the ordered sequence of eigenvalues of these matrices.*

Taking into account that the known mutual arrangement of CP's roots having certain multiplicities $q_\alpha$ uniquely determines this $q_\alpha$ by the correspondent eigenvalue's number in ordered sequence, one has to conclude that the next Statement is fulfilled:

**Statement 10.** *The class containing the orbit of the Hermitian matrix $P$ is determined by sets (56) constructed for all $\alpha = \overline{1, l}$.*

## Conclusion

We have proposed method that allows to split the set of orbits of Hermitian matrices into equivalence classes. The classes are labeled by a set of multiplicities that correspond to the ordered sequence of eigenvalues of these matrices. We found an opportunity to describe these classes using rational functions of the coefficients of the characteristic polynomial of the matrix. The explicit expression for syzygies is found in the case when the eigenvalues of have multiplicities more than 1 and when some of these multiplicities are repeated.

When studying the problem we have factorized characteristic polynomial of the matrix collecting all eigenvalues having the same multiplicities in one factor – polynomial with coefficients rationally expressed through CP's coefficients. In addition, we have found (with any accuracy) minimal distance between eigenvalues, as well as minimal and maximal eigenvalues of the matrix, using only rational functions of coefficients of CP. Then we have splitted the region containing all eigenvalues of CP onto uniform cells of maximal length such that each cell contains at most one eigenvalue. We can check whether the certain cell is empty or not using only rational functions of CP's coefficients.

## Acknowledgment

We are grateful to Professor V. P. Gerdt, to A. Khvedelidze, S. Yevlakhov, G. Giorgadze and D. Karkashadze for helpful and fruitful discussions.